\documentclass{article}
\usepackage[]{inputenc}
\usepackage{amssymb}
\usepackage{amsmath}
\usepackage{amsthm}
\usepackage{authblk}
\usepackage[]{asymptote}
\pdfoutput=1
\usepackage[
backend=biber,
style=alphabetic,
sorting=ynt
]{biblatex}
\usepackage[]{geometry}
\title{Restricted Positional Games}
\author[1]{Pranav Avadhanam}%put affiliations in - your high school etc[I just added my high school for affiliation - PA]
\author[2, 3]{Siddhartha G. Jena}
\affil[1]{Monta Vista High School}
\affil[2]{Princeton University Department of Molecular Biology}
\affil[3]{Harvard University Department of Stem Cell and Regenerative Biology}
\affil[ ]{\textit {\{pavadhanam3.1@gmail.com, sjena@fas.harvard.edu\}}}
\newtheorem*{definition}{Definition}
\newtheorem*{lemma}{Lemma}
\newtheorem*{theorem}{Theorem}
\begin{document}
\newtheorem*{corollary}{Corollary}
\maketitle
\section{Introduction}
A (two-player) \emph{strong positional game} played on a hypergraph $(X,H)$ ($H$ being a collection of subsets of $X$) is a game where two players sequentially label(``claim'') %the quotation marks are correct here, make sure they are correctly oriented throughout. also, consider italicizing all definitions throughout. -SJ [I believe I fixed all the quotation issues, and tried italicizing most of the definitions - PA]
vertices of $X$ with colors assigned to each player until all vertices of the board are claimed(see \cite{PosG} for a thorough introduction to positional games). The first player to claim all elements of a winning set $h \in H$ wins. If no player claims all the elements of a winning set, then the game results in a draw. Each particular sequence of vertices claimed by either player is called a \emph{play}, and the \emph{terminal position} is the final coloring(we use the terms \emph{coloring} and \emph{labeling} interchangeably) of the board as a result of a given play. \emph{Partial plays} are prefixes of plays(a partial play of length $k$ is the first $k$ moves of a given play), and the partial labeling that is the result of a partial play is called a \emph{position}. One such example of a positional game is \emph{Tic-Tac-Toe}, where the board is the $3$-by-$3$ grid: 
\begin{figure}[h]
    \centering
    \includegraphics{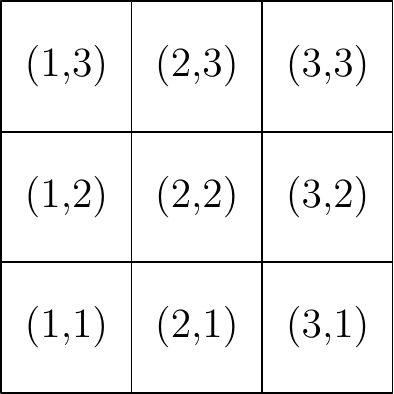}
    \label{fig:my_label}
\end{figure}
\\ Here, the collection of winning sets, $H$, is the set of rows, columns, and diagonals on the board. The first player in Tic-Tac-Toe labels their claimed elements by $``X"$, and the second player labels $``O"$. Denote this game as the $3^2$-game and define $[n]=\{i \in \mathbf{N}: 1\leq i \leq n\}$. The \emph{$n^d$ Tic-Tac-Toe} game is played on the $[n]^d$ hypercube with winning sets being \emph{geometric lines}. Given a board $[n]^d$, a geometric line in the board, $g$, is a collection of vertices $\{v_1, v_2, ..., v_n\}$, such that there exists an $n \times d$ matrix where each row vector is a distinct element of $g$, and where each column vector is an arithmetic progression(with common difference $-1$, $0$, or $1$). Positional games such as Tic-Tac-Toe are \emph{perfect information}, meaning that each player has knowledge of all the previous moves up to the current point in the play. Given that the set of available vertices to claim on each move is known, it is possible to simulate every possible play(searching through all branches of the game tree) and computationally determine optimal strategies for both players. The problem with this approach, when applied to games such as $n^d$ Tic-Tac-Toe, is that there are $(n^d)!$ possible plays to search through, and so the problem of finding (winning) strategies becomes intractable for small values of $n$ and very small values of $d$. This phenomenon, known as \emph{combinatorial explosion}(sometimes called ``combinatorial chaos"), is the primary motivator for analyzing these games combinatorially rather than computationally. While Tic-Tac-Toe has already been extensively studied from this viewpoint,%consider putting in a transition sentence here -SJ [For a transition, I added some writing about combinatorial chaos - PA]
\ a popular variant of Tic-Tac-Toe known as \emph{Connect-4}(that is played on a $6$-by-$7$ grid oriented vertically) provides an important example of what might be called a ``restricted positional game": games in which the set of available vertices on each move is somehow restricted due to additional constraints such as gravity. The winning sets in Connect-4 are any $4$ vertices of the board that are aligned in a row, column, or diagonal. Although the two games are defined similarly, Connect-Tac-Toe additionally has the influence of gravity to affect what plays are possible(with chips sliding down the columns to the available vertex of lowest height). Connect-4, like Tic-Tac-Toe, is a solved game(see \cite{C4Thesis}). %this is the result that you show, right? the way you have put it here makes it sound like a previous result! -SJ [I added my reference - PA]
We similarly generalize Connect-4 to the $n^d$ game played on the $[n]^d$ hypercube with winning sets being geometric lines, and the set of plays being restricted in a way similar to how gravity slides chips to the bottom of a column in ``real-world" play. 
\section{Connect-Tac-Toe}
\begin{definition}[Connect-Tac-Toe]
Connect-Tac-Toe is a strong positional game. It is played on a board $X$, $X = [n]^d$, for a given side-length $n$ and dimension $d$, and with a family of winning sets $W$, where $W$ is the set of all geometric lines in $X$. For each fixed $(x_1, x_2, x_3, ..., x_{d-1})$ in $[n]^{d-1}$, 
let $$ \textbf{c}_{x_1x_2...x_{d-1}} = \{(x_1, x_2, ..., x_{d-1}, k): k\in [n]\},$$ one of $n^{d-1}$ generalized 'columns' in $d$ dimensions. The game-play is as follows: for each partial play of $m$ moves, the ${(m+1)}^{th}$ move must be selected from a set of available vertices $A$. Let $P$ be the set of previously claimed vertices. Then, 
\newline
$$\mathcal{A} = \{\mathbf{a}: \exists \mathbf{c}_{x_1, x_2, ..., x_{d-1}}, \mathbf{a} \in \mathbf{c}_{x_1, x_2, ..., x_{d-1}} \land \mathbf{a}_d = \min\{c_d: c \in \mathbf{c}_{x_1, x_2, ..., x_{d-1}}\cap (X \backslash P) \}\}\footnote{For a given vector $\mathbf{x}$, we denote ${
\mathbf{x}}_l$ as the $l^{th}$ coordinate of $\mathbf{x}$. Similarly, $\mathbf{x}_{i,j}$ is the $j^{th}$ coordinate of vector $\mathbf{x}_i$},$$ $\mathcal{A}$ being the set of unclaimed vertices with the least vertical coordinate amongst the unclaimed vertices in their column. 
\end{definition}
Given the following partial play $((1,1,1),(1,1,2),(2,2,1),(3,3,1),(1,1,3))$ in the $3^3$-game, we mark the set of available vertices to be claimed with $A$(column $\textbf{c}_{1,1}$ is marked in green):
\begin{figure}[h]
    \centering
    \includegraphics{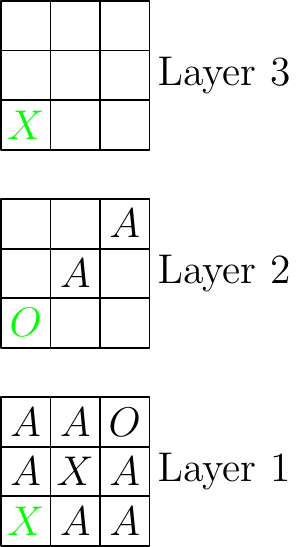}
    \label{fig:my_label}
\end{figure}
\section{Restrictions on Positional Games}
$C2T\footnote{We will have $C2T$ denote the game Connect-Tac-Toe, and $3T$ denote the game Tic-Tac-Toe.}$ is an example of a strong positional game where for a given partial play the set of available moves for the next turn(there are $\vert X \vert$ turns total) is restricted past just the set of unclaimed vertices in the board. A ``usual" strong positional game might be represented by the game hypergraph $(X, F)$, where $X$ is the board and $F$ is the collection of winning sets in $X$. Let $ \mathcal{P}_i $ be the set of partial plays of length $i$, where $1 \leq i \leq \vert X \vert$, then $\mathcal{P} = \bigcup_{i=1}^{\vert X \vert}{\mathcal{P}_i}$. A \emph{restricted positional game} is instead represented by the 3-tuple $(X, F, \mathcal{A})$, where the additional function $\mathcal{A}$ maps from the set of all partial(but potentially full) plays $\mathcal{P}$ to the board $X$. $\mathcal{A}$ is interpreted as assigning the set of available moves given a partial play up to a certain point in the game(much like the previous definition of $\mathcal{A}$ in section $2$) and  we say that $\mathcal{A}(\emptyset) = X$. In ``usual" positional games $\mathcal{A}(m_1, m_2, ... m_k)$ is assumed to equal $X \backslash \{m_1, m_2, ... m_k\}$, but in $C2T$ $\mathcal{A}$ has an alternative definition and where $\mathcal{A}(m_1, m_2, ... m_k)$ is some subset of $X\backslash \{m_1, m_2, ... m_k\}$. The sole difference between $3T$ and $C2T$ is the definition of $\mathcal{A}$, and therefore $C2T$ can be seen as a certain restricted form of $3T$. In particular, letting $Pl_{3T}(n,d)$(resp. $Pl_{C2T}(n,d)$) be the set of plays for $n^d$ $3T$(resp. $C2T$), we note that ${Pl_{C2T}(n,d)} \subseteq Pl_{3T}(n,d)$. Note that $Pl_{3T}(n,d)=(n^d)!$, since each play is an ordering on the vertices of the board. To make it clearer that $C2T$ is indeed a restricted form of $3T$, we have the following lemma:
\begin{lemma}
$\vert Pl_{C2T}(n,d)\vert = \frac{Pl_{3T}(n,d)}{(n!)^{n^{d-1}}}$
\end{lemma}
\begin{proof}
Note that for a given partial play $p$, $\mathcal{A}(p)$ contains at most $1$ vertex from each column. $\mathcal{A}(p)$ has no vertices in a column if and only if that column has all of its vertices already claimed(that column has been ``chosen" $n$ times previously). Because choosing a not completely claimed column will uniquely determine the vertex that will be selected for that move, there is a bijection between $Pl_{C2T}(n,d)$ and the set of permutations of the multiset $M$ with $i$(each $i$ representing a column) of multiplicity $n$(eventually choosing the column $n$ times in the course of the play) for every $i \in [n^{d-1}]$ and $M$ containing no other elements. As a classical result, $\vert M \vert = \binom{n^d}{n,n,...n}$($n^{d-1}$ $n$'s in the bottom argument). So we have $\vert Pl_{C2T}(n,d) \vert = \vert M \vert = \frac{(n^d)!}{(n!)^{n^{d-1}}}$.
\end{proof}
If we define $TP_{3T}(n,d)$ to be the set of terminal positions in $[n]^d$, we note that $\vert TP_{3T}(n,d) \vert = {n^d \choose \lfloor \frac{n^d}{2} \rfloor}$(elements of $TP_{3T}(n,d)$ are known as \emph{halving colorings} of $[n]^d$). Consider the \emph{layers} $L_i$, $L_i=\{(x_1,x_2,...x_{d-1}, i)\in [n]^d\}$, of the $n^d$ hypercube. Any halving coloring with $L_1$ monochromatic in color $``O"$ is necessarily not in $TP_{C2T}(n,d)$(defined similarly to $TP_{3T}(n,d)$) since the first vertex which would be marked $``X"$ must be in $L_1$ if $C2T$-restrictions were in place. So we know that $\vert TP_{C2T}(n,d) \vert<\vert TP_{3T}(n,d)\vert$ for $n,d>1$. Here is an example of a terminal position that is in $TP_{3T}(4,2)$ but not in $TP_{C2T}(4,2)$:
\begin{figure}[h]
    \centering
    \includegraphics{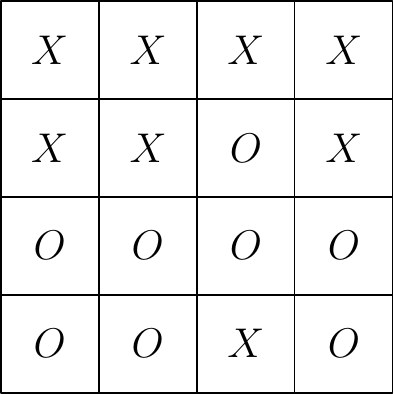}
    \label{fig:my_label}
\end{figure}
\section{The $C2T$ Hales Jewett Number}
The \emph{Hales-Jewett Theorem}, proved in \cite{RegPos}, provides a condition where as long as the dimension is sufficiently large($\geq HJ(n)$), there will be no drawing position in the Tic-Tac-Toe game on $[n]^d$. Pairing this with the following theorem, we have what Beck, Pegden, and Vijay \cite{HJ} have called a ``soft existential criterion'': without any demonstrated strategy, the properties of the game hypergraph alone guarantee a first player win. As a way to distinguish $TP$ and $C2T$ game-theoretically, we study variants of this number. Alternatively, and perhaps more to the point, we can consider the win number $w(n)=\min(d_0:\forall d\geq d_0, \ first \ player's \ win \ on \ n^d\ 3T)$ and $w_{C2T}(n)=\min(d_0:\forall d\geq d_0, \ first \ player's \ win \ on \ n^d\ C2T)$, but these seem to be more difficult to study.
\begin{theorem}[Strategy-Stealing]
The first player always has a drawing or winning strategy for $n^d$ Tic-Tac-Toe.
\end{theorem}
The reasoning for this comes from the classical ``Strategy Stealing'' argument, due to Nash. Suppose the second player had a winning strategy. Then we can show that the first player can have a first move at random without decreasing the chance of the first player winning or drawing. After the first player having ``given up'' the first move, the second player becomes the effective ``first player'' and the first player can adopt whatever winning strategy the second player had to win the game. Since this implies the second player can never have a winning strategy, we have the result.
\begin{theorem}[Hales-Jewett]
For every positive integer $n$, there exists a minimum integer $HJ(n)$, such that for every $2$-coloring of $[n]^{HJ(n)}$, there exists a monochromatic geometric line.
\end{theorem}
Call a coloring \emph{proper} if it does not contain a monochromatic geometric line. All though the Hales-Jewett theorem shows that there exists a dimension at which all $2$-colorings of $[n]^d$ are improper, the Hales-Jewett theorem implies something much stronger to hold for $HJ(n)$.
\begin{corollary}[$HJ(n)$ is a strict threshold]
$\forall d \geq HJ(n)$, every $2$-coloring of $[n]^d$ is improper and $\forall d < HJ(n)$, there exists a proper $2$-coloring of $[n]^d$.
\end{corollary}
\begin{proof}
By definition, for every $d<HJ(n)$, there exists a proper $2$-coloring of $[n]^d$. Suppose for some dimension $d\prime$, we have that every $2$-coloring of $[n]^{d\prime-1}$ is improper. Given that $[n]^{d\prime}$ consists of $n$ layers that are ``copies" of $[n]^{d\prime-1}$, any coloring of $[n]^{d\prime}$ with a layer necessarily containing a monochromatic geometric line will also contain a monochromatic geometric line as a whole coloring: so every $2$-coloring of $[n]^{d\prime}$ is improper. Setting $d\prime=HJ(n)+1$, we verify that $\forall d \geq HJ(n),\ \nexists \ proper \ $2$-coloring$ of $[n]^{d}$ by induction on $d\prime$.
\end{proof}
Define $HJ_{1/2}(n)=\min(d:\forall \ c\in TP_{3T}(n,d), \exists \ monochromatic \ geometric \ line)$ and $HJ_{C2T}(n)=\min(d:\forall \ c\in TP_{C2T}(n,d), \exists \ monochromatic \ geometric \ line)$. Note that $HJ(n)\geq HJ_{1/2}(n) \geq HJ_{C2T}(n)$, this following from $TP_{C2T}(n,d) \subset TP_{3T}(n,d)$. For these variants of the Hales-Jewett Number($HJ(n)$) such a ``threshold" property has not been proven. In the paper \cite{HJ}(which also discussed the threshold property), exponential lower bounds were proven for the Hales-Jewett and Hales-Jewett halving number:
\begin{theorem}[Beck-Pegden-Vijay]
$HJ(n)\geq HJ_{1/2}(n) \geq \frac{2^{(n-2)/4}}{3(n-2)^4}$
\end{theorem}
We are able to provide logarithmic lower bounds for the Hales-Jewett $C2T$ number, and it is an open question whether these bounds can be improved to exponential. If these logarithmic bounds are tight, we can quantitatively describe just how restrictive(from more of a game-theoretic perspective) $C2T$ is versus $3T$. For a given $2$-coloring $C:X\rightarrow \{0,1\}$, let $C\prime:X\rightarrow \{0,1\}$ denote the \emph{color-flip} of $C$ with $C:x \mapsto  0$ if and only if $C\prime:x \mapsto 1$.
\begin{lemma}
Every $2$-coloring of $[n]^d$ such that $\lceil\frac{n}{2}\rceil$ layers are colored by a halving coloring $C$ and $\lfloor\frac{n}{2}\rfloor$ layers are colored by $C\prime$ is an element of $TP_{C2T}(n,d)$.
\end{lemma}
\begin{proof}
Let $C$ be halving coloring of the set of $n^{d-1}$ columns in $[n]^d$(alternatively $C$ represents a halving coloring of a layer of the hypercube). Denote our assignment of colorings to each layer by $f:[n]\mapsto\{C,C\prime\}$. If $n$ is even, then $C\prime$ is also a halving coloring. Call a (partial) play $C2T$-valid, if it satisfies $C2T$-restrictions. Then we may construct a $C2T$-valid play layer-by-layer: let $C=(\{x_1,x_2,...x_{\frac{n^{d-1}}{2}}\},\{y_1,y_2,...y_{\frac{n^{d-1}}{2}}\})$. For even $n$, we have the $C2T$-valid play(here we have presented it as a sequence of vertices rather than columns) on $n^d$, $P=(p_1,p_2,...p_{n^d})$, such that $\forall i\in\{0,1,...n-1\}$ with $f(i)=C$(resp. $f(i)=C\prime$), $(p_{n^{d-1}\cdot i+1},p_{n^{d-1}\cdot i+2},...p_{n^{d-1}\cdot i+n^{d-1}})=(x_{1,i},y_{1,i},...x_{\frac{n^{d-1}}{2},i},y_{\frac{n^{d-1}}{2},i})$ (resp. $(y_{1,i},x_{1,i},...y_{\frac{n^{d-1}}{2},i},x_{\frac{n^{d-1}}{2},i})$). For odd $n$, we claim that $\forall w \in [n]$, $\exists$ $C2T$-valid partial play on $\bigcup_{i=1}^{w}L_i \backslash X$ where $X$ is a monochromatic(with respect to $f(w)$) subset of $L_w$ such that $\vert\bigcup_{i=1}^{w}L_i\backslash X\vert$ is even and having a position labeled as $(f(1),...f(w-1),f(w)[L_w\backslash X])$. In this case, let $C=(\{x_1,x_2,...x_{\lceil\frac{n^{d-1}}{2}\rceil}\},\{y_1,y_2,...y_{\lfloor\frac{n^{d-1}}{2}\rfloor}\})$. We proceed with induction on $w$. For the base case of $w=1$, we explicitly construct a $C2T$-valid partial play $P_1$ on $L_1\backslash \{\phi_{1}\}$ where for $f(1)=C$(resp. $f(1)=C\prime$) $\phi \in$ $\{x_1,x_2,...x_{\lceil\frac{n^{d-1}}{2}\rceil}\}$(resp. $\{y_1,y_2,...y_{\lceil\frac{n^{d-1}}{2}\rceil}\}$): we may arbitrarily choose $\phi$ to be $x_{\lceil\frac{n^{d-1}}{2}\rceil}$(resp. $y_{\lceil\frac{n^{d-1}}{2}\rceil}$) so that $P_1 = (x_1, y_1,...,x_{\lfloor\frac{n^{d-1}}{2}\rfloor},y_{\lfloor\frac{n^{d-1}}{2}\rfloor})$(resp. $(y_1, x_1,...,y_{\lfloor\frac{n^{d-1}}{2}\rfloor},x_{\lfloor\frac{n^{d-1}}{2}\rfloor})$). Suppose for some $w-1$, we have that there is a $C2T$-valid partial play $P_{w-1}$ on $\bigcup_{i=1}^{w-1}L_i \backslash X$(with an appropriate $X$ and position of $P_{w-1}$ being a subset of $f$). Note that for the terminal position $f(1),f(2),...f(w-1)$, the difference in sizes between the color classes is $\vert\vert\{1\leq i \leq w-1: f(i)=C\}\vert-\vert\{1\leq i \leq w-1: f(i)=C\prime\}\vert\vert$. Since the terminal position of $\bigcup_{i=1}^{w-1}L_i \backslash X$ has zero difference in the size of the color classes, $X$ must be responsible for this overall discrepancy: $\vert X \vert = \vert\vert\{1\leq i \leq w-1: f(i)=C\}\vert-\vert\{1\leq i \leq w-1: f(i)=C\prime\}\vert\vert$. Given that for $f$ there are at most $\lceil \frac{n}{2} \rceil$ layers that can be $C$(and at most $\lfloor \frac{n}{2} \rfloor$ that can be $C\prime$), $\vert X \vert \leq \lceil \frac{n}{2} \rceil$. There are at least $\lfloor\frac{n^{d-1}}{2}\rfloor-\lceil\frac{n}{2}\rceil>\vert X \vert$ columns $\mathbf{c}$ with $\mathbf{c}_w$ of the opposite color and $X \cap \mathbf{c} = \emptyset$ such that we can choose $\mathbf{c}_w$ to pair with an element from $X$. We want to extend the partial play $P_{w-1}$ to a partial play $P_{w}$ on $\bigcup_{i=1}^{w}L_i\backslash X\prime$. For each of the elements  in $X$, we have shown that we can pair them with an above element from $L_w$ while satisfying the $C2T$-restrictions. But after pairing off(sequentially selecting vertices with opposite labels in the terminal position) these elements of $X$, there will remain at least(depending on the color of $X$ and what $f(w)$ is) $\lfloor\frac{n^{d-1}}{2}\rfloor - \lceil\frac{n}{2}\rceil$ vertices of each color on $L_w$ that have not yet been claimed in the course of the partial play. Pairing off these remaining vertices, we will be left with a monochromatic set $X\prime \subset L_w$ and $P_{w}$ will be $C2T$-valid on $\bigcup_{i=1}^{w}L_i\backslash X\prime$. Having completed the induction, we know the case for $w=n$ holds: there exists a $C2T$-valid play $P$ on $[n]^d \backslash X_n$. Furthermore, we know that $\vert X_n \vert = \vert\vert f^{-1}[C]\vert-\vert f^{-1}[C\prime]\vert\vert = \vert f^{-1}[C]\vert-\vert f^{-1}[C\prime]\vert = 1$, so $X_n$ must contain a single element $\chi$ of $\{x_{1,n},x_{2,n},...x_{\lceil\frac{n^{d-1}}{2}\rceil,n}\}$(resp. $(y_{1,n},y_{2,n},...y_{\lceil\frac{n^{d-1}}{2}\rceil,n})$) for $f(n)=C$(resp. $f(n)=C\prime$). Selecting $\chi$ as the first player's final move of the play, we have that $(P,col(\chi))$$\footnote{Let $col(x)$ be the column a vertex $x\in[n]^d$ is in.}$ is a $C2T$-valid play(represented as a sequence of columns) on $[n]^d$ with terminal position described by $f$.
\end{proof}
\begin{lemma}[Folklore]
There are $\frac{(n+2)^d-n^d}{2}$ geometric lines in the $[n]^d$ board.
\end{lemma}
\begin{proof}
Consider the geometric line $g=\{v_1,v_2,...v_n\}$. There exists two orderings $(v_1,v_2,...v_n),(v_n,v_{n-1},...v_1)$ on $g$ that we count such that for every $i \in [n]$, $c_i=(v_{1,i},v_{2,i},...v_{n,i})$ is equal to $(1,2,...n)$, $(n,n-1,...1)$, or $(a,a,...a)$(for some $a \in [n]$). Given that there are $n$ values of $a$, and two additional choices for whether $(c_i)$ is an arithmetic progression of positive or negative common difference, there are $(n+2)^{d}$ (possibly degenerate) orderings on the geometric lines. Subtracting for the $n^d$ degenerate geometric lines(lines where all the columns are arithmetic progressions with common difference 0), and dividing by $2$ to take into account the fact that there are twice as many orderings than geometric lines, there must be $\frac{(n+2)^d-n^d}{2}$ unique geometric lines.
\end{proof}
\begin{theorem}
We have the following lower bound:
$$HJ_{C2T}(n) > \frac{n-6}{4\log_2(n)}$$
\end{theorem}
\begin{proof}
$dim(L_i)=d-1<HJ_{C2T}(n)$, so there exists some proper halving coloring, $C$, of each $L_i$. Note that it suffices to show(given the previous lemma) that there exists a proper coloring of $n^d$ with $\lceil\frac{n}{2}\rceil$ layers colored $C$ and $\lfloor\frac{n}{2}\rfloor$ layers colored $C\prime$ provided some restriction on $n$. Since $C$ is proper, $C\prime$ is also proper, and therefore all geometric lines contained within a layer will be non-monochromatic. Suppose $L_i$ is colored with $C$ and $L_j$ is colored with $C\prime$. Then, for every $x_i\in L_i$, the corresponding $x_j\in L_j$ that is in the same column as $x_i$ will have coloring opposite of that of $x_i$. Therefore the column containing both $x_i$ and $x_j$ will be guaranteed to be non-monochromatic. Spanning over all $x_i$, all columns will also be guaranteed to be non-monochromatic no matter $C$. We know there are $\frac{(n+2)^d-n^d}{2}$ geometric lines. The total number of geometric lines that we still haven't accounted(to be guaranteed non-monochromatic) is going to be $\frac{(n+2)^d-n^d}{2}-(n\cdot\frac{(n+2)^{d-1}-n^{d-1}}{2})-n^{d-1}=(n+2)^{d-1}-n^{d-1}$. All of these remaining geometric lines $G_1, G_2, ... G_k(k=(n+2)^{d-1}-n^{d-1})$ intersect every layer once. Fix $L_1$ to be colored by $\mu$(denoted by $(L_1,\mu)$, $\mu \in \{C,C\prime\}$). Let $proj_{i}(j, \ell)$ be the unique $g \in L_i$ such that $g$ and $G_j(\ell)$($\ell$th vertex of the geometric line) are in the same column. Let $x_1=\vert\{j \in [k]: \mu(proj_{1}(j,1))=\mu(proj_{1}(j,2))\}\vert$, and $y_1 = k-x_1$. If we choose $(L_2,\mu)$, the set of geometric lines that are currently monochromatic in our partially constructed coloring is reduced by $y_1$. If we choose $(L_2,M)$($M$ being the color-flip of $\mu$), $k$ is instead reduced by $x_1$. Denote the remaining quantity of monochromatic geometric lines as $k\prime$. We may force $k\prime=k-max(x_1,y_1)\leq\frac{k}{2}$, reducing $k$ by a factor of at least 2. We repeat this process, coloring $L_{i+1}$ in terms of $L_i$, until layer q, where $k^{(q-1)}=0$. Since we only have freedom to choose colorings between $C,C\prime$ in $\lfloor\frac{n}{2}\rfloor$ layers(we cannot arbitrarily color between $C,C\prime$ in $\lceil\frac{n}{2}\rceil$ layers since there is a possibility that all of those $\lceil\frac{n}{2}\rceil$ layers would be colored with $C\prime$ which would violate our initial assumption), $q=\log_2(k=(n+2)^{d-1}-n^{d-1})+2\leq \lfloor\frac{n}{2}\rfloor$. If $d=HJ_{C2T}(n)$, then there does not exist a proper coloring in $TP_{C2T}(n,d)$, and so we have $\log_2((n+2)^{HJ_{C2T}(n)-1}-n^{HJ_{C2T}(n)-1})+2 > \lfloor\frac{n}{2}\rfloor$. From here, we derive the lower bounds.
\end{proof}
\printbibliography
\end{document}